\documentclass{article}
\usepackage{amsmath,enumerate,graphicx}
\usepackage{url}
\usepackage{subcaption}
\usepackage{float}
\usepackage{xcolor}

\newtheorem{theorem}{Theorem}[section]

\newtheorem{definition}[theorem]{Definition}

\newtheorem{remark}[theorem]{Remark}

\numberwithin{equation}{section}

\begin{document}
%
\title{Singular points of envelopes and offsets:\\ a transition from 2D to 3D}
%
%
\maketitle  

\begin{center}
\author{Thierry Dana-Picard$^{1,2}$, Arie Haenel$^{1}$, Daniel Tsirkin$^{1}$\\
$^1$ Jerusalem College of Technology, Havaad Haleumi st. 21\\ Jerusalem 9116011, Israel\\
$^2$ Jerusalem Michlala College, Duvdevani St. 24\\ Jerusalem, Israel\\
ndp@jct.ac.il,arie.haenel@g.jct.ac.il,academictsirkin@gmail.com}
\end{center}

\small
\begin{center}
The authors' names are listed in alphabetical order. \\ All three authors have
contributed at the same level to this work.
\end{center}
\normalsize

\begin{abstract}
We recall properties of the offsets of a (plane) kiss curve
and their singularities, then we make a transition to 3D, considering
the surface obtained by revolution of the kiss curve about one of its
axes. When possible, geometric constructions using automated methods
are performed. This induces a need to focus on algebraic computations.
Networking between different kinds of software is of utmost importance,
in order to use their respective strengths (dynamic geometry, dragging
points, strong algebraic computations, etc.). A core difference between
the 2D and 3D situations is that automated commands existing in 2D
are not available in 3D. This makes the DGS less relevant in 3D than
in 2D. In 3D, most of the work must be performed with the CAS. The
singular points of a kiss curve in the plane are tacnodes (cusps) and
crunodes (self-intersection). The singularities are not isolated anymore
and cuspidal edges appear. In this setting, implicitization is a heavy
endeavour and we need another technology, related to AI; we describe
this briefly. The 2D activities have been proposed to in-service teachers
learning towards an advanced degree.
\end{abstract}

\section{Introduction}

\subsection{Parametric presentations and cartesian equation of a surface}

This paper reports on an ongoing work about various purely geometric constructions which yield plane algebraic curves, studied in a technology-rich environment. Here we consider offsets of a plane curve\footnote{Leibnitz called them \emph{parallel curves} \cite{leibnitz}. Nowadays, this name may be misleading, as shown in certain descriptions by generative AI; see \cite{kiss curve}.} and generalize the study to offsets of an algebraic surface in 3-dimensional space. The curve and the surface (called the progenitors) are initially defined by parametric presentations involving trigonometric functions. Using the substitution  
\begin{equation}
\label{substitution trig to rational}
\begin{cases}
\cos t =\frac{1-u^2}{1+u^2}\\
\sin t = \frac {2u}{1+u^2}
\end{cases}
\end{equation}
followed eventually  by algebraic manipulations, this can be transformed into a rational parametrization. If needed, this can be be transformed into a polynomial equation of the surface. This polynomial equation can be obtained either via the computation of an elimination ideal \cite{cox} or in one step using Maple's \textbf{implicitize} command. More details in next subsection.

\subsection{Automated methods}
Automated methods are developed and implemented in software \cite{ThEdu}, requesting both Dynamic Geometry (for interactive exploration and conjectures) and a Computer Algebra System (for algebraic manipulations). We use GeoGebra\footnote{Freely downloadable from \url{http://www.geogebra.org}}, and  its companion package GeoGebra-Discovery\footnote{Freely downloadable from \url{https://github.com/kovzol/geogebra-discovery}}. The 2nd one has symbolic implementations of commands, which work numerically in the standard GeoGebra (such as \textbf{Locus}, \textbf{Envelope}).  GeoGebra has an embedded CAS, called Giac \cite{kovacs parisse}, but when more computational power is needed, we use Maple 2024. These automated commands have no counterpart for 3D Geometry, therefore we have to focus on algebraic computations. Maple also has automated commands for implicitization, parametrization and looking for singularities of a curve, but all these are devoted to plane curves only.

Moreover, it is well-known that software has hard time to plot a curve in a neighborhood of a singular point. Here, GeoGebra-Discovery command \textbf{Plot2D} is an efficient tool to overcome this pitfall. 

Both parametric presentations and implicit equations are used, switching from parametric to polynomial representation being based on packages for Gr\"obner bases and Elimination. We wish to emphasize the importance of the networking between different kinds of software, and wish to see more possibilities of communication between them. Since \cite{eugenio}, a few advances have been performed, but we have still to copy-paste by hand between GeoGebra and Maple. In Section \ref{implicitization - networking between technologies}, we present an alternative implicitization method using Sylvester resultants within a Python/SymPy framework, illustrating cross-platform collaboration.

What we present switches mathematics education from the traditional sequence 
\begin{equation*}
definition \; \rightarrow \; theorem  \; \rightarrow \; proof  \; \rightarrow \; examples
\end{equation*}
towards a more experimental approach: 
  \begin{equation*}
exploration \; \rightarrow \; conjecture \; \rightarrow \; proof \; \rightarrow \; theorem
 \end{equation*}
 
 This is quite clear for geometry, but can be adapted to other domains. We experienced this with a group of in-service mathematics teachers learning towards a 2nd degree in Mathematics Education. At the end of the course, they expressed their happiness to have discovered a new approach to mathematics.

\section{Envelopes and offsets in the plane}
Envelopes of parameterized families of plane curves are important, both for the mathematics involved and for applications \cite{safety}. Offsets of a plane curve $\mathcal{C}$ are sometimes considered as a special case of envelopes, but this is not accurate. Due to Euler, offsets are also called parallel curves, a misleading name. Envelopes and offsets are generally different objects, similar only in particular cases, e.g., when the curve  $\mathcal{C}$ is smooth and the offset distance small enough. Studying their topology may be fascinating, as the topology of an offset can be more complicated than that of its progenitor \cite{alcazar-sendra}.

\begin{definition}
\label{def offset of a plane curve}
Let be given a plane curve $\mathcal{C}$ and a positive real number. At every non singular point $A$ of  $\mathcal{C}$, we define a unit normal vector $\vec{N_A}$. The offset of $\mathcal{C}$ at distance $d$ is the geometric locus of points $M$ such that $\overrightarrow{AM}=d \cdot \vec{N_A}$.
\end{definition}

\subsection{Offsets of an ellipse}

\begin{figure}[htb]
\begin{center}
\includegraphics[width=3.5cm]{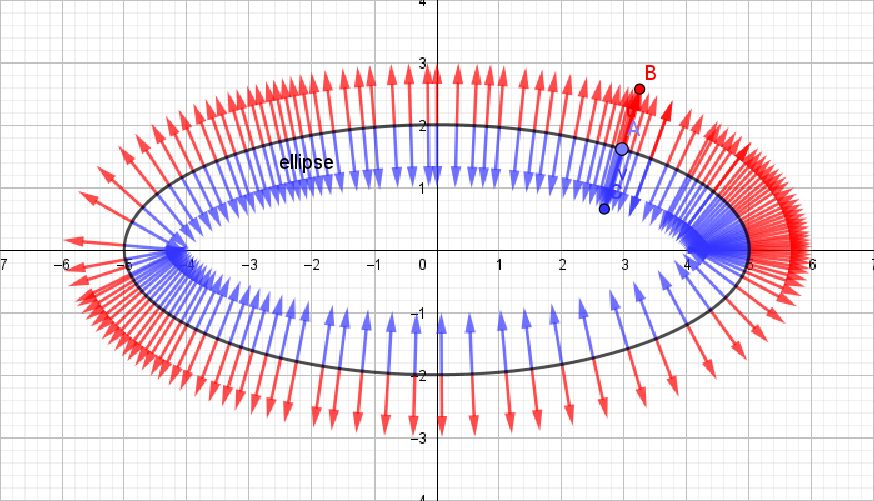}
\quad
\includegraphics[width=3.5cm]{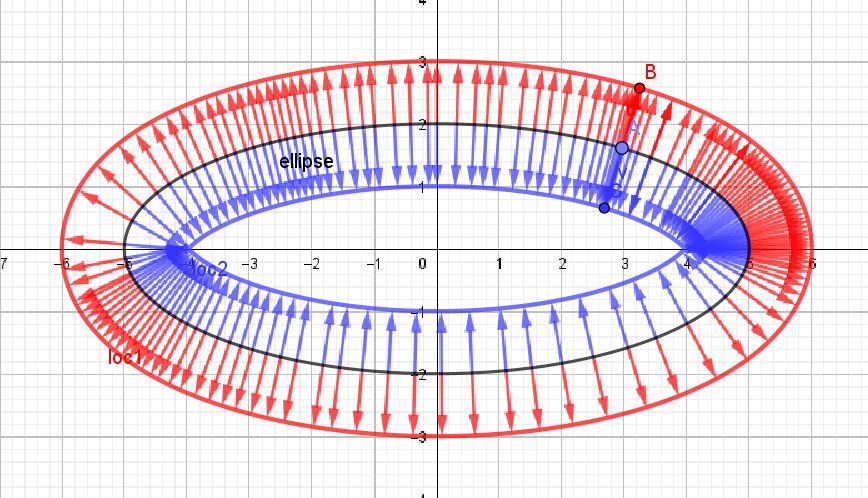}
\quad
\includegraphics[width=3.5cm]{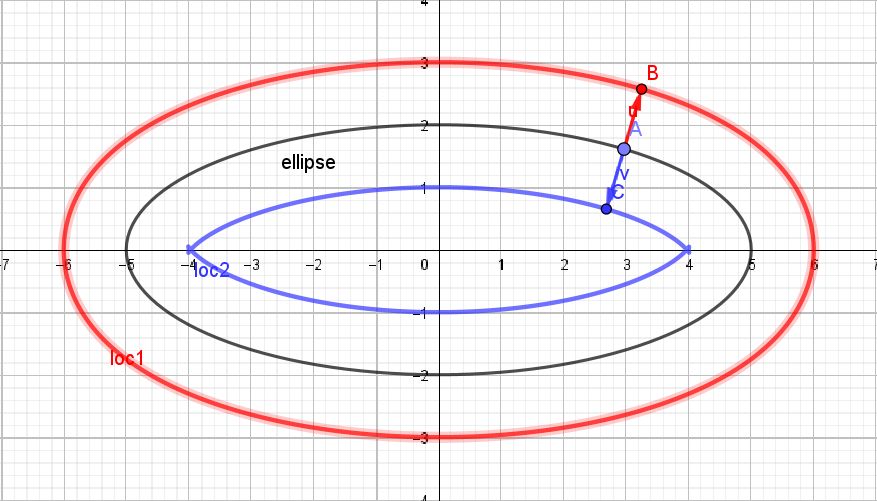}
\caption{The two components of an offset of an ellipse}
\label{fig 1}
\end{center}
\end{figure}

We consider an ellipse whose equation is $\frac{x^2}{25}+\frac{y^2}{4}=1$. The exploration of this subsection is made with a GeoGebra applet \url{https://www.geogebra.org/m/ngc4unxk}. An interactive construction and exploration of the offset at distance 1 of this ellipse is illustrated in Figure \ref{fig 1}. From left to right:
\begin{enumerate}[(a)] 
\item First trial: move the point $A$ along the ellipse.  Some vectors, their number being determined by the speed of the hand on the mouse, and obtained with \emph{Trace On}, emphasize the dynamic aspect of the work.
\item Second step: the command \textbf{Locus(Point,Point)} has been added. Its output consists in a curve, plotted without providing an equation. Such an equation cannot be obtained with GD, as the construction involves steps which are not supported by the \textbf{LocusEquation} command. A popup gives the reason of the command reject.
\item Here we display the offset only. It is the union of two disjoint loops, one can be called outer offset, the other one inner offset. Note that it seems that something "interesting" occurs with the inner component.
\end{enumerate}  

\begin{figure}[htb]
\begin{center}
\includegraphics[width=5cm]{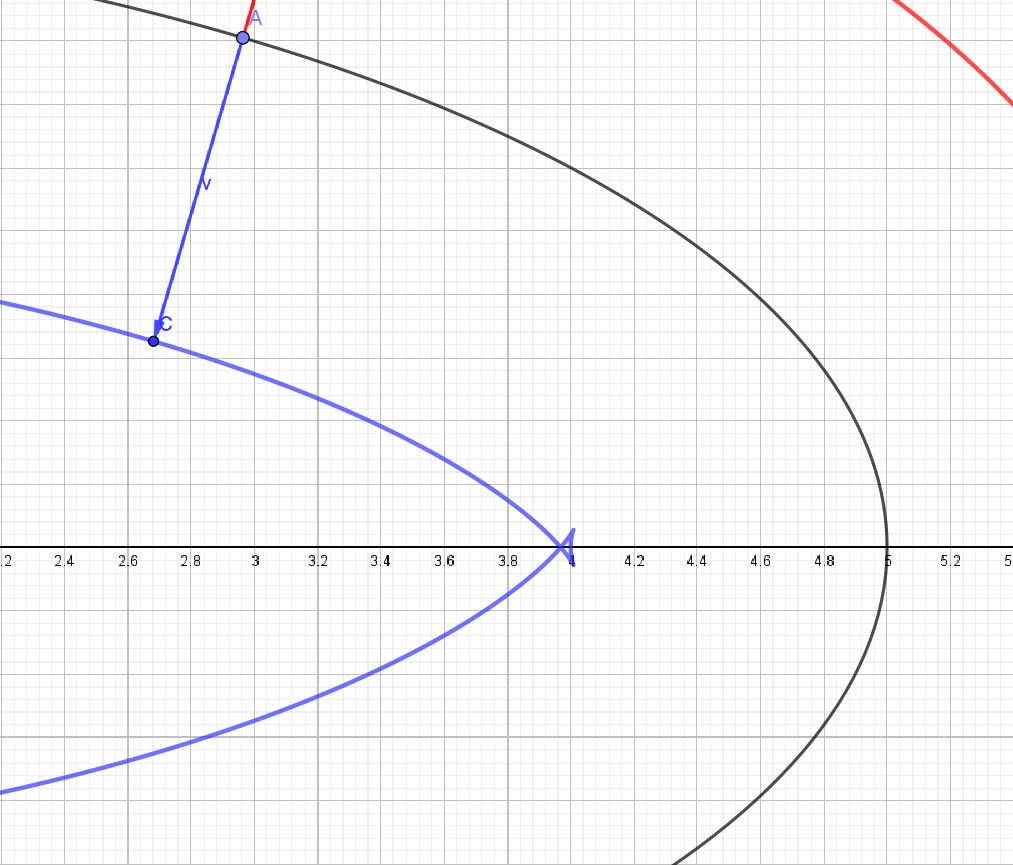}
\caption{The inner component has self-intersection and cusps}
\label{fig 1 zoomed}
\end{center}
\end{figure}

Strong zooming may reveal what happens with the inner component, as shown in Figure \ref{fig 1 zoomed}. In particular, cusps (also called tacnodes) appear.

\begin{figure}[htb]
\begin{center}
\includegraphics[width=3.5cm]{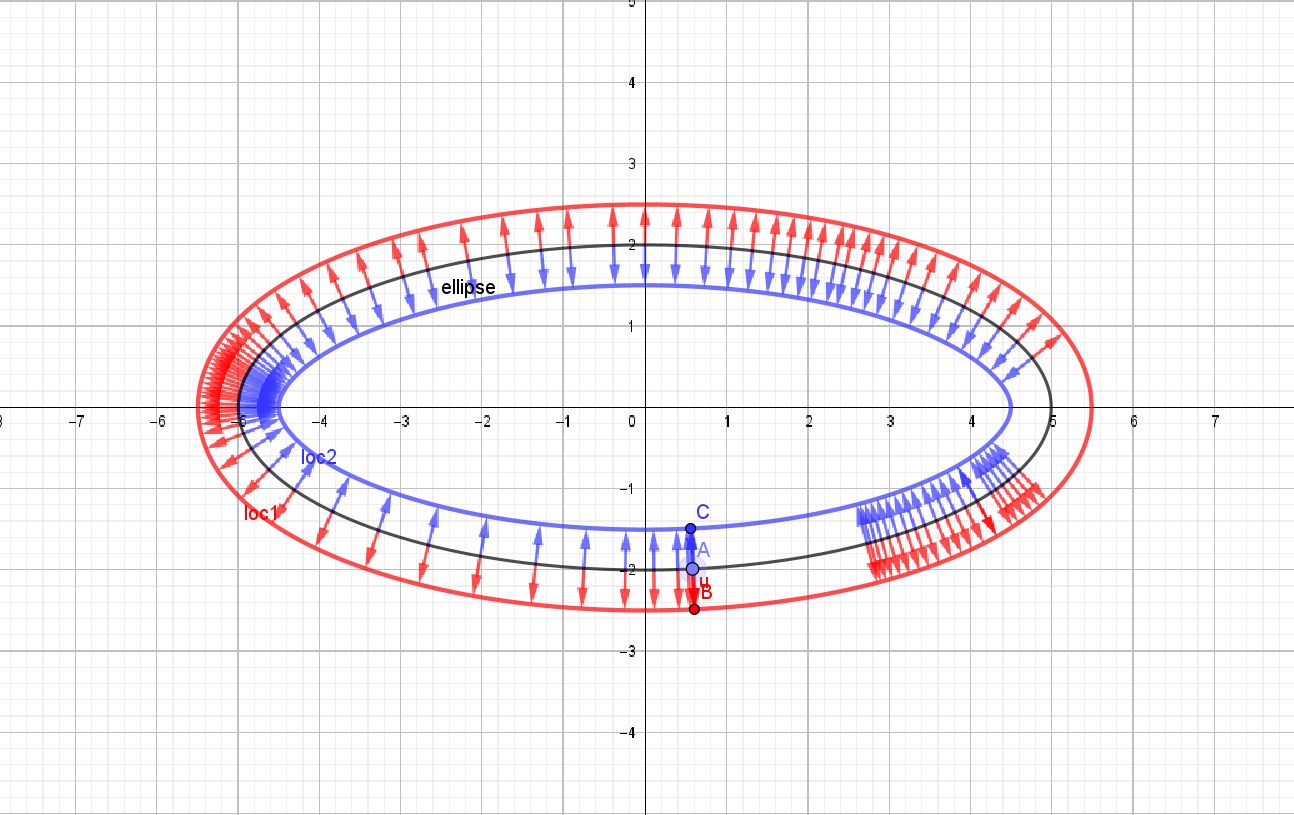}
\quad
\includegraphics[width=3.5cm]{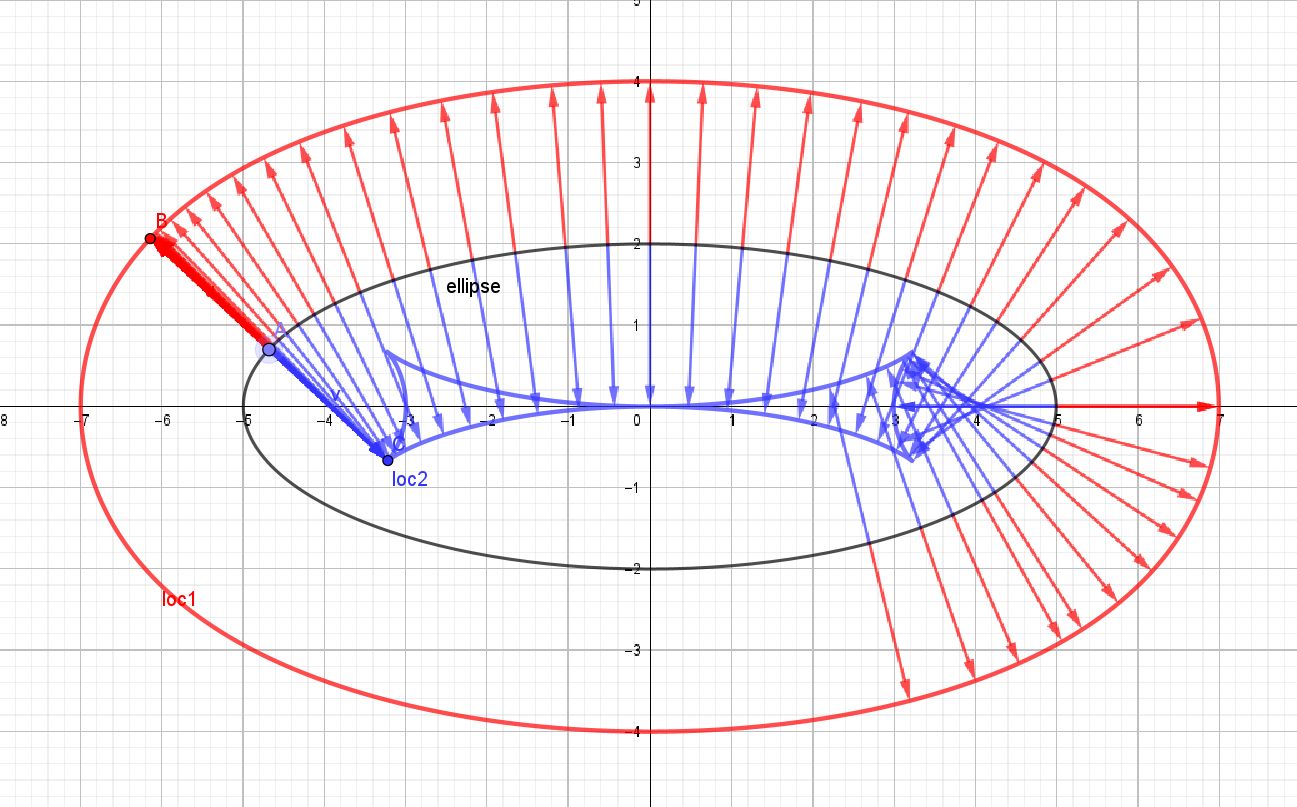}
\quad
\includegraphics[width=3.5cm]{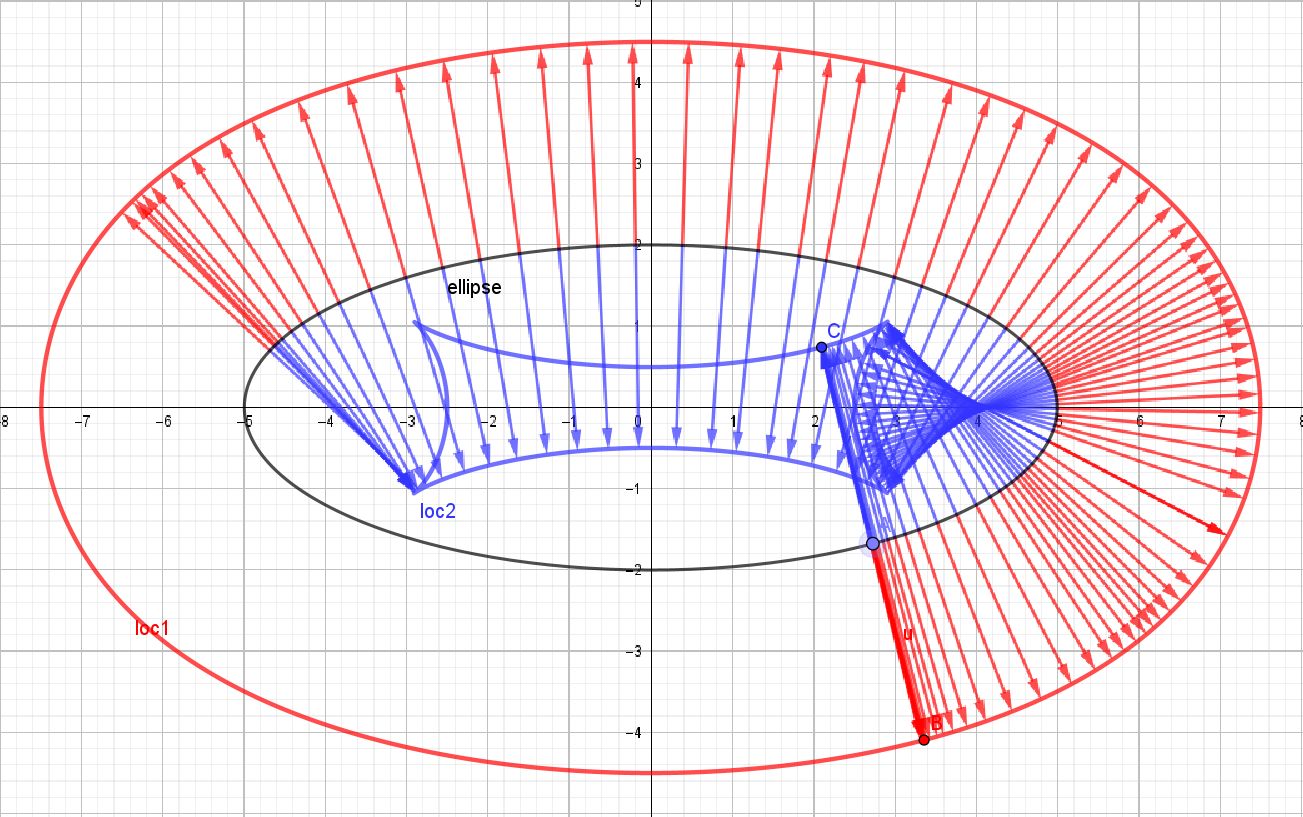}
\caption{The topology of the offset depends on the distance}
\label{fig 1 dependence on the distance}
\end{center}
\end{figure}

The dynamic exploration reveals that the existence of tacnodes (cusps)  depends on the offset distance, as illustrated in Figure \ref{fig 1 dependence on the distance}. To obtain it, we used a segment to parameterise the distance, as generally GD requires a pure geometric construction. Some of the commands still work when the data is provided by an equation, but not always. From left to right, the distance is equal to $0.5$, $2$ and $2.5$ respectively.

For $d=0.5$, the offset seems to be smooth. For $d=1$, crunodes and cusps appear on the inner offset (v.s. Figure \ref{fig 1}). For $d=2$, the inner offset exhibits cusps on the side, and seems to be self-tangent at the origin. For $d=2.5$, the inner components looks as a curved rectangle, with 4 cusps and no crunode.

\begin{figure}[htb]
\begin{center}
\includegraphics[width=3.5cm]{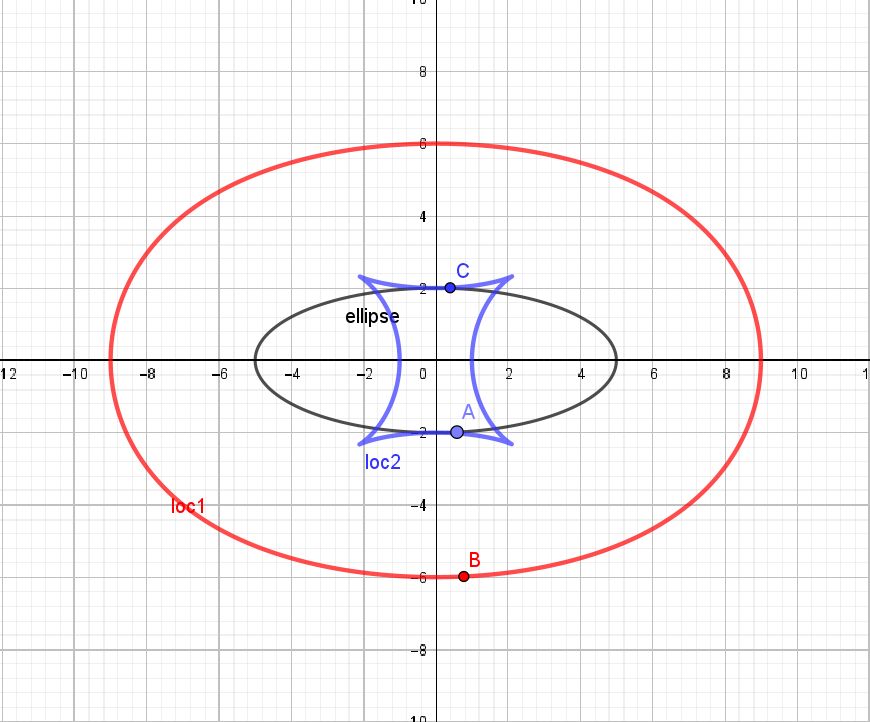}
\quad
\includegraphics[width=3.5cm]{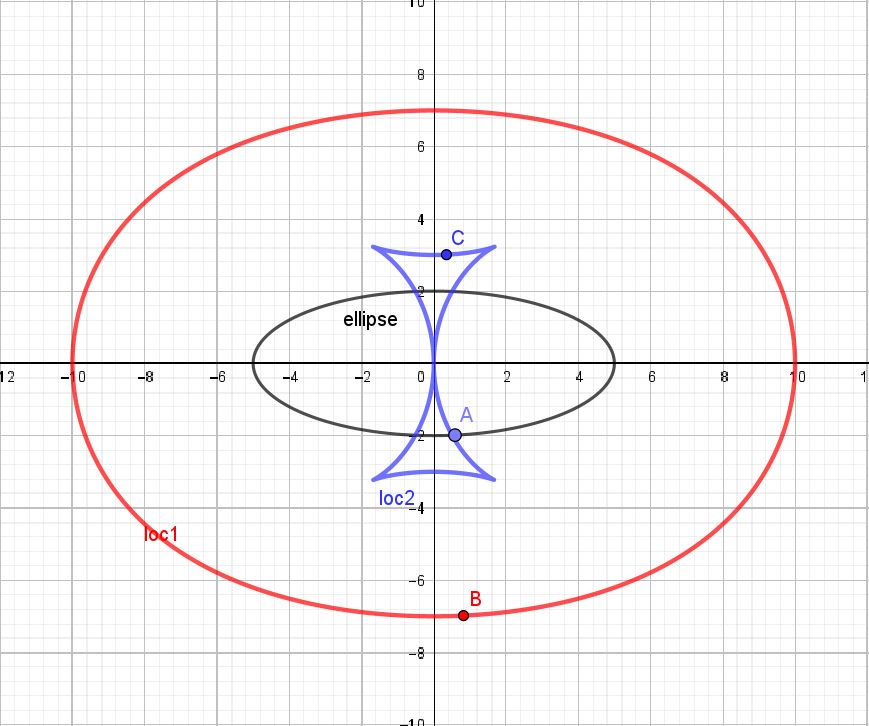}
\quad
\includegraphics[width=3.5cm]{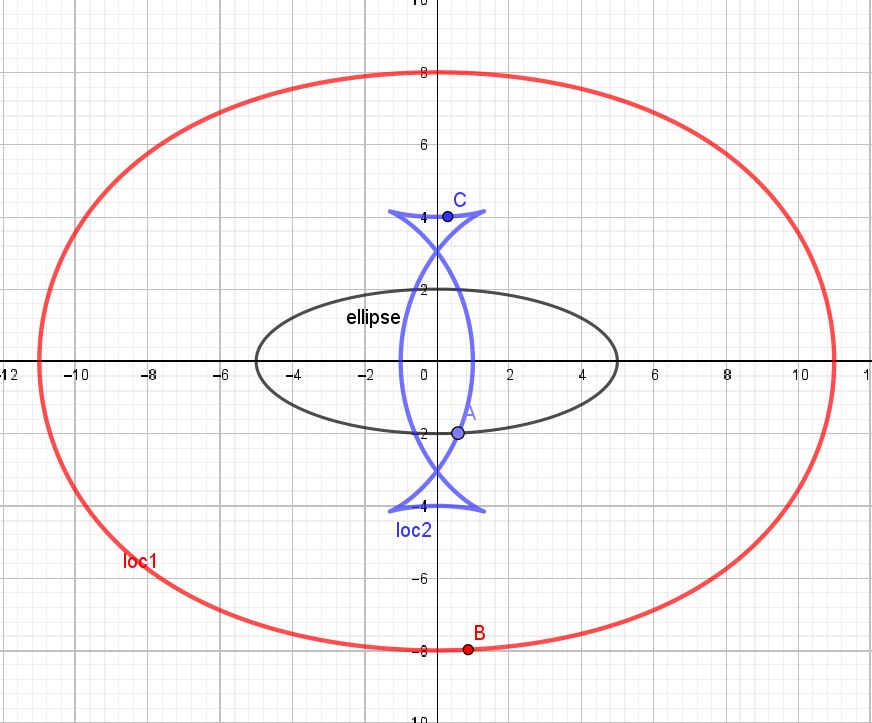}
\caption{The topology of the offset depends on the distance}
\label{fig 1 larger distances}
\end{center}
\end{figure}

If the offset distance increases, the situation can change more. Figure \ref{fig 1 larger distances} illustrates this for $d=4,5,6$ (from left to right). Note that the main axis of the inner offset is now vertical, and also that the name "inner offset" is not accurate anymore.

All the observations above have to be confirmed by algebraic computations. In particular, the limiting values of the distance should be determined. This will request applications of some of the methods with Maple, and already described in \cite{kiss curve}.

\subsection{Offsets of a kiss curve}
\begin{figure}
\begin{center}
\includegraphics[width=3cm]{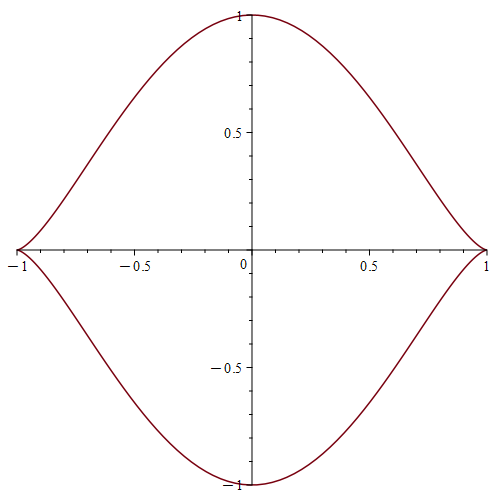}
\qquad \qquad  
\includegraphics[width=3.5cm]{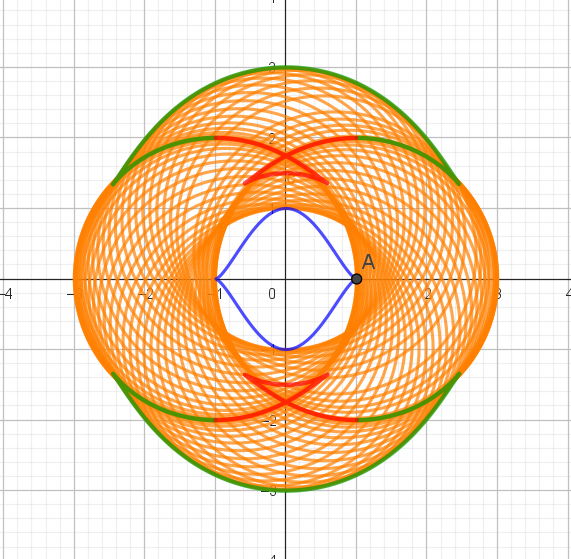}
\caption{A GeoGebra plot to compare the envelope of unit circles and the offset at distance 1}
\label{fig offset_env_comp}
\end{center}
\end{figure}

A kiss curve is plane curve $\mathcal{C}$ given by the following parametric presentation:
\begin{equation}
\label{kiss_curve_eq}
\begin{cases}
\mathcal{C}_x (t) = \cos t\\
\mathcal{C}_y (t) =  a \sin^3 t 
\end{cases}
\end{equation}
where $a $ is a positive real parameter. WLOG, we work with $a=1$, the topology of the curve is the same for every value of $a$. 
Envelopes of circles of constant radius centered on $\mathcal{C}$ and offsets of $\mathcal{C}$ have been described and their singular points characterized in  \cite{kiss curve}, together with Maple code for the algebraic computations. Figure \ref{fig offset_env_comp} shows a kiss curve and a screenshot of a GeoGebra session to compare the envelope of the family of unit circles centred on $\mathcal{C}$ and the offset at distance 1.

\begin{figure}[htb]
\begin{center}
\includegraphics[width=3cm]{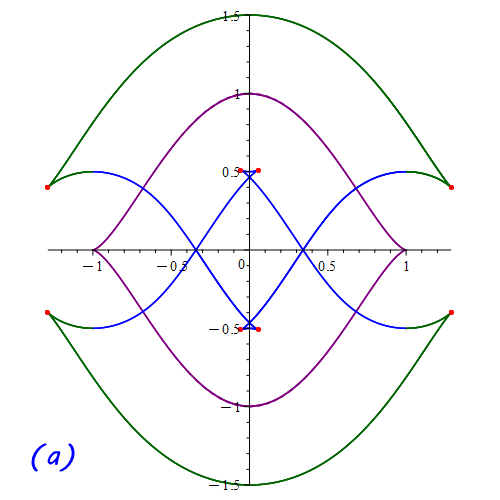}
\qquad
\includegraphics[width=3cm]{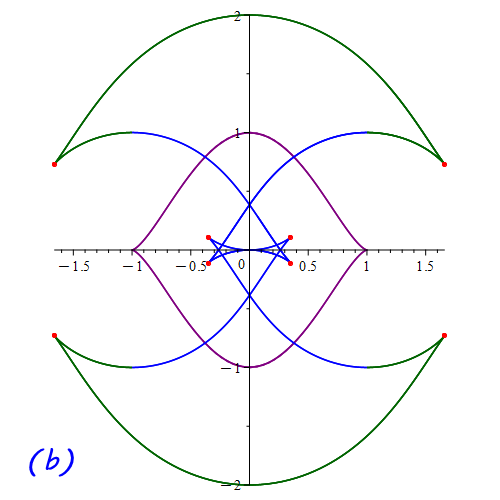}
\qquad
\includegraphics[width=3cm]{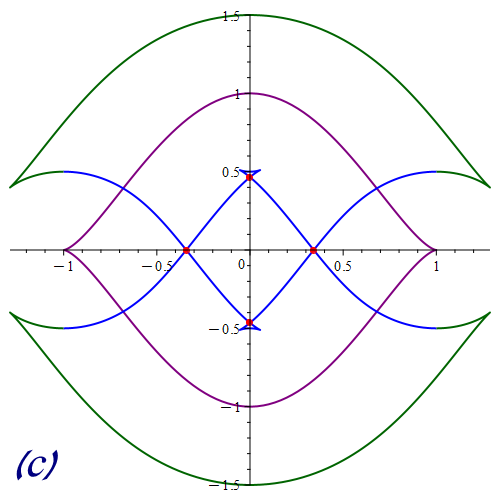}
\caption{Singular points of an offset}
\label{fig 2}
\end{center}
\end{figure}

The outer and inner offsets at distance $d$ is defined by the following parameterizations:

\begin{equation}
\label{offsets_parametrization}
\begin{aligned}
\text{Offset}_{\text{in}}(t) &=
\begin{cases}
x(t) = \mathcal{C}_{x}(t) - d \cdot \dfrac{d}{dt}\mathcal{C}_{y}(t),\\
y(t) = \mathcal{C}_{y}(t) + d \cdot \dfrac{d}{dt}\mathcal{C}_{x}(t),
\end{cases}
\\[1em]
\text{Offset}_{\text{out}}(t) &=
\begin{cases}
x(t) = \mathcal{C}_{x}(t) + d \cdot \dfrac{d}{dt}\mathcal{C}_{y}(t),\\
y(t) = \mathcal{C}_{y}(t) - d \cdot \dfrac{d}{dt}\mathcal{C}_{x}(t).
\end{cases}
\end{aligned}
\end{equation}
Two approaches can be used to determine the singular points of the offset:
\begin{enumerate}
\item Find the points, for which  $\frac{dx}{dt} = \frac{dy}{dt} = 0$. 
\item Find the points, for which the offset curvature is not defined.
\end{enumerate}
The 1st approach determines all the candidates to be singular, the 2nd approach determines cusps, not crunodes (points of self-intersection). To find the crunodes, the following system of equations has to be solved:

\begin{equation}
\begin{cases}
\text{Offset}_{\text{in}}(t) =\text{Offset}_{\text{in}}(s)\\
\text{Offset}_{\text{out}}(t) =\text{Offset}_{\text{out}}(s)
\end{cases}
\end{equation}
Note that these equations take into account the fact that the outer and inner offsets do not intersect, but at the points where they are "glued". In Figure \ref{fig 2}, the singular points are enhanced (in red).

\section{Transition to 3D: a virtual clam}
We consider now a generalization to 3D of the previous study. The main issue is that most of the singularities will not be isolated anymore. 

\subsection{A virtual clam}
\begin{figure}[htb]
\begin{center}
\includegraphics[width=6cm]{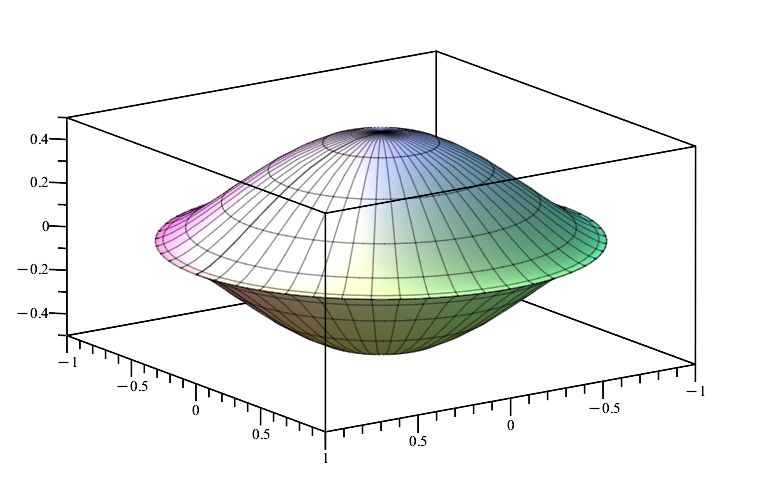}
\caption{A virtual clam}
\label{fig clam param}
\end{center}
\end{figure}

We consider the surface (we call it a clam) given by the parametric presentation
\begin{equation}
\label{param eq clam}
\begin{cases}
x= \cos u \cos v\\
y=  \cos u \sin v\\
z= \frac 12 \sin^3 u
\end{cases}
\end{equation}
The coefficient 1/2 has been chosen for aesthetic reasons. The topology of the surface is the same for every positive value of the parameter $a$: the transformation changing the surface for a value $a_1$ to the surface for a value $a_2$ is an affinity, i.e. a homeomorphism. The surface for $a=1/2$ is displayed in Figure \ref{fig clam param}.

To derive a cartesian (symbolic) equation for the surface, we run the following Maple code. The core of the process is the transformation of the trigonometric expressions into rational expressions, using the same substitution as in \cite{safety}. 
\footnotesize
\begin{verbatim}
restart;
with(plots);
with(PolynomialIdeals);
setoptions(scaling = constrained);
xclam := cos(u)*cos(v);
yclam := cos(u)*sin(v);
zclam := 1/2*sin(u)^3;
clamplot := plot3d([xclam, yclam, zclam], u = 0 .. 2*Pi, v = 0 .. Pi,
 scaling = constrained);
xc := subs(cos(u)=(-s^2 + 1)/(s^2 + 1), subs(cos(v)=(-t^2 + 1)/(t^2 + 1),
           xclam));
yc := subs(cos(u)=(-s^2 + 1)/(s^2 + 1), subs(sin(v)=2*t/(t^2 + 1),yclam));
zc := subs(sin(u) = 2*s/(s^2 + 1), zclam);
P1 := denom(xc)*x - numer(xc);
P2 := denom(yc)*y - numer(yc);
P3 := denom(zc)*z - numer(zc);
J := <P1, P2, P3>;
JE := EliminationIdeal(J, {x, y, z});
F := Generators(JE)[1];
evala(AFactor(F));
Fclam := x^6 + 3*x^4*y^2 + 3*x^2*y^4 + y^6 - 3*x^4 - 6*x^2*y^2 - 3*y^4 
         + 3*x^2 + 3*y^2 + 4*z^2 - 1:
implicitplot3d(Fclam = 0, x = -1 .. 1, y = -1 .. 1, z = -1 .. 1,
               scaling = constrained, numpoints = 20000, 
               style = patchnogrid)
\end{verbatim}
\normalsize

Actually, elimination provides a polynomial of degree 8, which can be factored into a polynomial of degree 6 and two linear factors. The two linear factors, namely $z-1/2$  and $z+1/2$ are irrelevant, as it appears immediately from plotting the octic surface and comparing with Figure \ref{fig clam param}. They appear because of topological reason, as explained in \cite{safety}.

\begin{remark}
The polynomial of degree 6 can be obtained also by a single command, namely \emph{implicitize}, in Maple's \emph{algcurves} package. This command requests a "guess" of the degree. If the guess is too small, no output is given. If it is too large, the computation time  may be much longer than for the accurate value and its output is a list  of polynomials, all of them multiple of the accurate one.
\end{remark}
From now on, we work with the following cartesian equation for the clam:
\begin{equation}
\label{cartesian eq clam}
x^6 + 3x^4y^2 + 3x^2y^4 + y^6 - 3x^4 - 6x^2y^2 - 3y^4 + 3x^2 + 3y^2 + 4z^2 - 1=0
\end{equation}
and denote the polynomial by $F(x,y,z)$.

\subsection{The family of unit spheres centred on the clam}
\begin{figure}[htb]
\centering
\begin{subfigure}{3.5cm}
 \includegraphics[width=3.5cm]{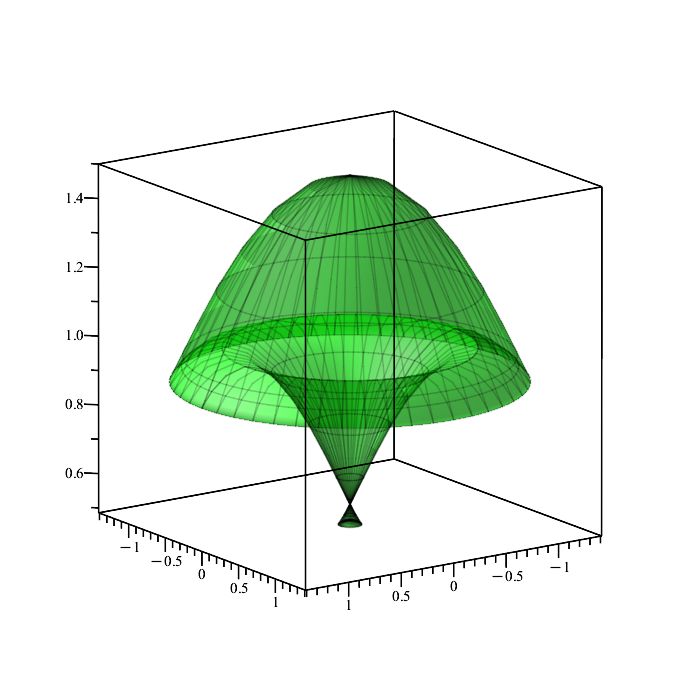}
\caption{from above} 
\end{subfigure}
\qquad
  \begin{subfigure}{3.5cm}
  \includegraphics[width=3.5cm]{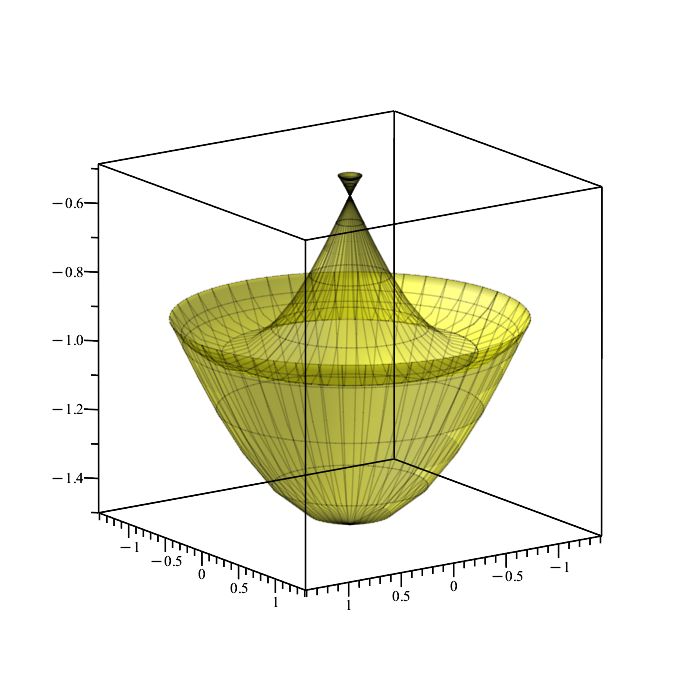}
 \caption{from below} 
\end{subfigure}
\qquad
\begin{subfigure}{3.5cm}
 \includegraphics[width=3.5cm]{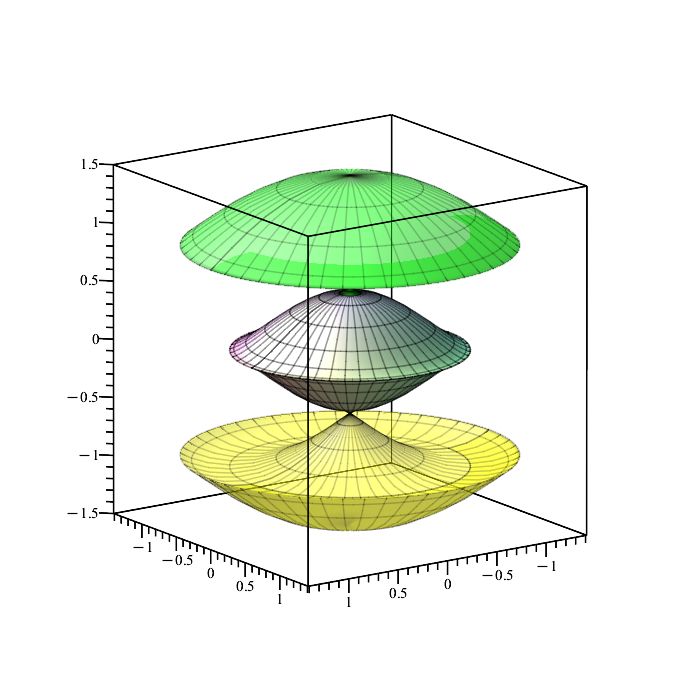}
 \caption{Clam and envelope}
 \end{subfigure}
 \caption{Envelope of unit spheres on a clam}
\label{fig env spheres on clam}
\end{figure}

We consider the family of unit spheres $\mathcal{S}_{u,v}$ centred on the clam. The spheres are defined by the equation $F(x,y,z,u,v)=0$, where:
\begin{equation}
\label{eq spheres}
F(x,y,z,u,v)=(x-\cos u \cos v)^2 +(y-  \cos u \sin v)^2+(z-\frac 12 \sin^3 u)^2-1.
\end{equation}
An envelope of the family $\{ \mathcal{S}_{u,v} \}$, if it exists, is defined by the solutions of the system of equations:
\begin{equation}
\label{eq system envelope}
\begin{cases}
F(x,y,z,u,v)=0\\
\frac{\partial F}{\partial u} (x,y,z,u,v)=0\\
\frac{\partial F}{\partial v} (x,y,z,u,v)=0\\
\end{cases}
\end{equation}
To make formulas more compact, we denote  $R = \sqrt{-9\sin^4(u) + 9\sin^2(u) + 4}$. The output gives parametric presentations of  2 surfaces:
\small
\begin{equation}
\label{param eq envelope spheres 1}
\begin{cases}
x = - \frac 14 \cos v \cos u \; \left( 3 \sin^4 u - \frac{3\sin u \; \left(R \; \sin^3 u - 4 \right)}{R + 4} - 4\right)\\
y = -\frac 14 \cos u \sin v \; \left( 3 \sin^4 u - \frac{3\sin u \; \left(R \; \sin^3 u - 4 \right)}{R + 4} - 4\right)\\
z = \frac{\sin^3 u (R - 4)}{2R}
\end{cases}
\end{equation}
\normalsize
and 
\small
\begin{equation}
\label{param eq envelope spheres 2}
\begin{cases}
x = - \frac 14 \cos v \cos u \; \left( 3 \sin^4 u - \frac{3\sin u \; \left(R \; \sin^3 u - 4 \right)}{R + 4} - 4\right)\\
y = -\frac 14 \cos u \sin v \; \left( 3 \sin^4 u - \frac{3\sin u \; \left(R \; \sin^3 u - 4 \right)}{R + 4} - 4\right)\\
z = \frac{\sin^3 u (R + 4)}{2R}
\end{cases}
\end{equation}
\normalsize
These surfaces are displayed separately in Figure \ref{fig env spheres on clam}(a) and (b), and together with the clam in  Figure \ref{fig env spheres on clam}(c).

\subsection{Implicitization}
\label{implicitization - networking between technologies}
As usual, the analytic definition of an envelope provided a parametric presentation of the envelope. It may be important to derive an implicit (cartesian) equation for the obtained surface. As the envelope is given by a trigonometric parametrization, it is well known that rational parametrizations exist, from which polynomial equations can be derived, using algorithms on Gr\"obner bases and elimination of the parameters. This has been done, for example, in \cite{hyperbolisms,AI-DPH}; a general description of an algorithm is given in \cite{lastra sendra sendra}.

The rational substitution approach leads to heavy computations with CAS in order to: a) derive polynomial equations, b) generate an ideal in a polynomial ring in 5 variables over the real field (2 of the variables being the parameters), and c) eliminate the parameters using Gr\"obner basis techniques via the \texttt{EliminationIdeal} command. 

We chose another approach with the Sylvester resultant method. This method performs implicitization of parametric envelope surfaces via a single Sylvester resultant eliminating the parameter $s = |\sin(u)|$ from two polynomial constraints in cylindrical coordinates $(\rho^2, z)$. The radical $R = \sqrt{-9\sin^4 u + 9\sin^2 u + 4}$ was eliminated algebraically by substituting $w = z - s^3/2$ and squaring the relation $w^2 = 4/R^2$, avoiding explicit radical manipulation. This method yields directly polynomials of degree 14  in $(\rho^2, z)$ for each envelope surface.\footnote{This computational approach and the implementation script were developed through interactive collaboration with an AI-assisted computational system (using Claude 4.5 with Model Context Protocol integration for SymPy), which proposed the key substitution strategy that simplified the parameter elimination process.}

The resulting implicit equations $F_9(\rho^2, z) = 0$ and $F_{10}(\rho^2, z) = 0$ contain 50 terms each and are too lengthy to reproduce in full. We show the leading terms:

\begin{align}
F_9(\rho^2, z) &= \frac{3^{20}}{2^{20}} (\rho^2)^{10} + \frac{3^{20}}{2^{18}} (\rho^2)^{9} z^{2} - \frac{3^{19} \cdot 109}{2^{21}} (\rho^2)^{9} \notag \\
&\quad + \frac{3^{21}}{2^{19}} (\rho^2)^{8} z^{4} - \frac{3^{19} \cdot 373}{2^{21}} (\rho^2)^{8} z^{2} + \cdots = 0 \label{eq:implicit9}\\
F_{10}(\rho^2, z) &= \frac{3^{20}}{2^{20}} (\rho^2)^{10} + \frac{3^{20}}{2^{18}} (\rho^2)^{9} z^{2} - \frac{3^{20} \cdot 31}{2^{21}} (\rho^2)^{9} \notag \\
&\quad + \frac{3^{21}}{2^{19}} (\rho^2)^{8} z^{4} - \frac{3^{20} \cdot 7 \cdot 17}{2^{21}} (\rho^2)^{8} z^{2} + \cdots = 0 \label{eq:implicit10}
\end{align}

where $\rho^2 = x^2 + y^2$.

\begin{figure}[htb]
\centering
\begin{subfigure}{5cm}
 \includegraphics[width=6cm]{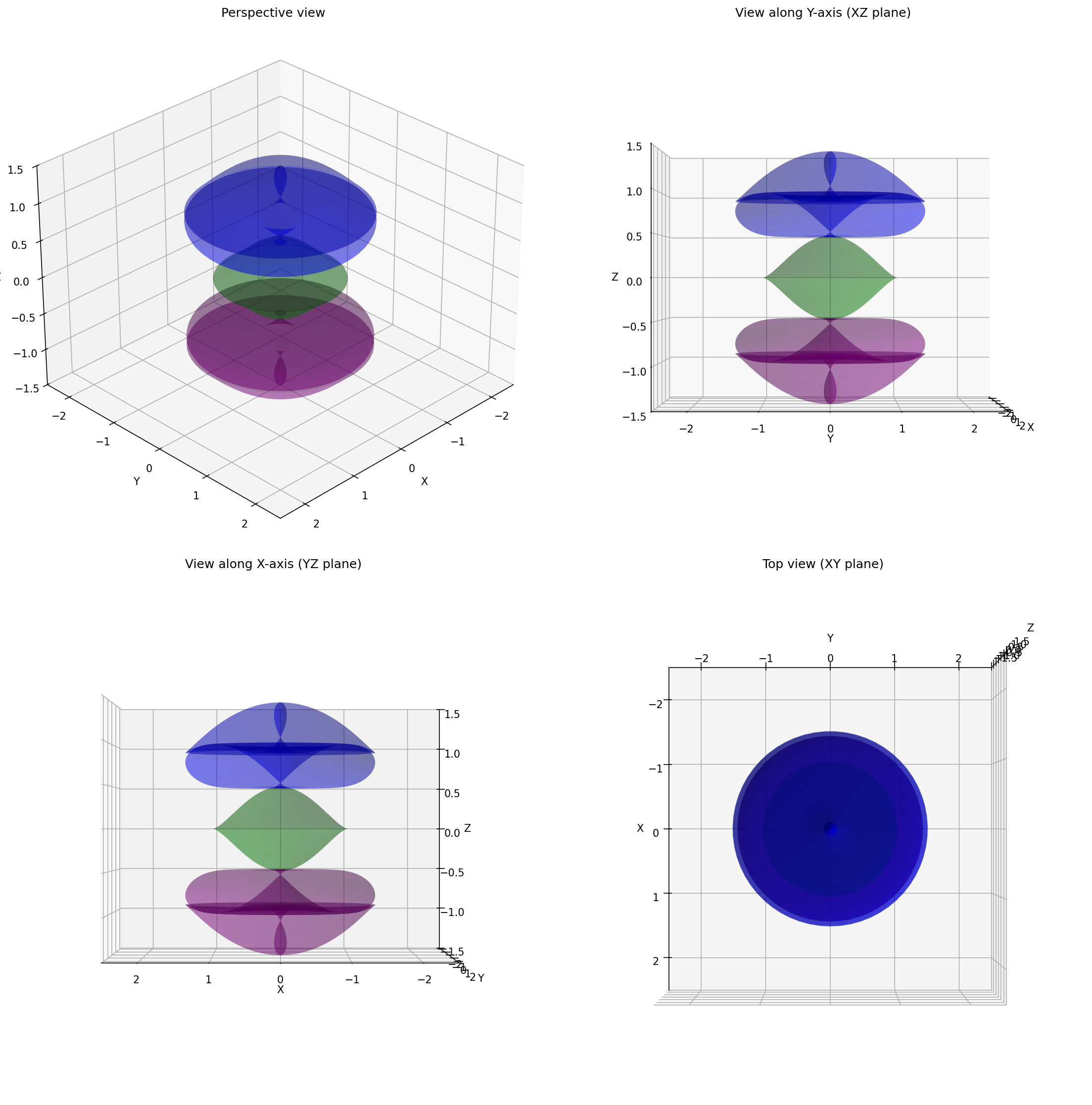}
\caption{Parametric Plot} 
\end{subfigure}
\qquad \qquad
  \begin{subfigure}{5cm}
  \includegraphics[width=6cm]{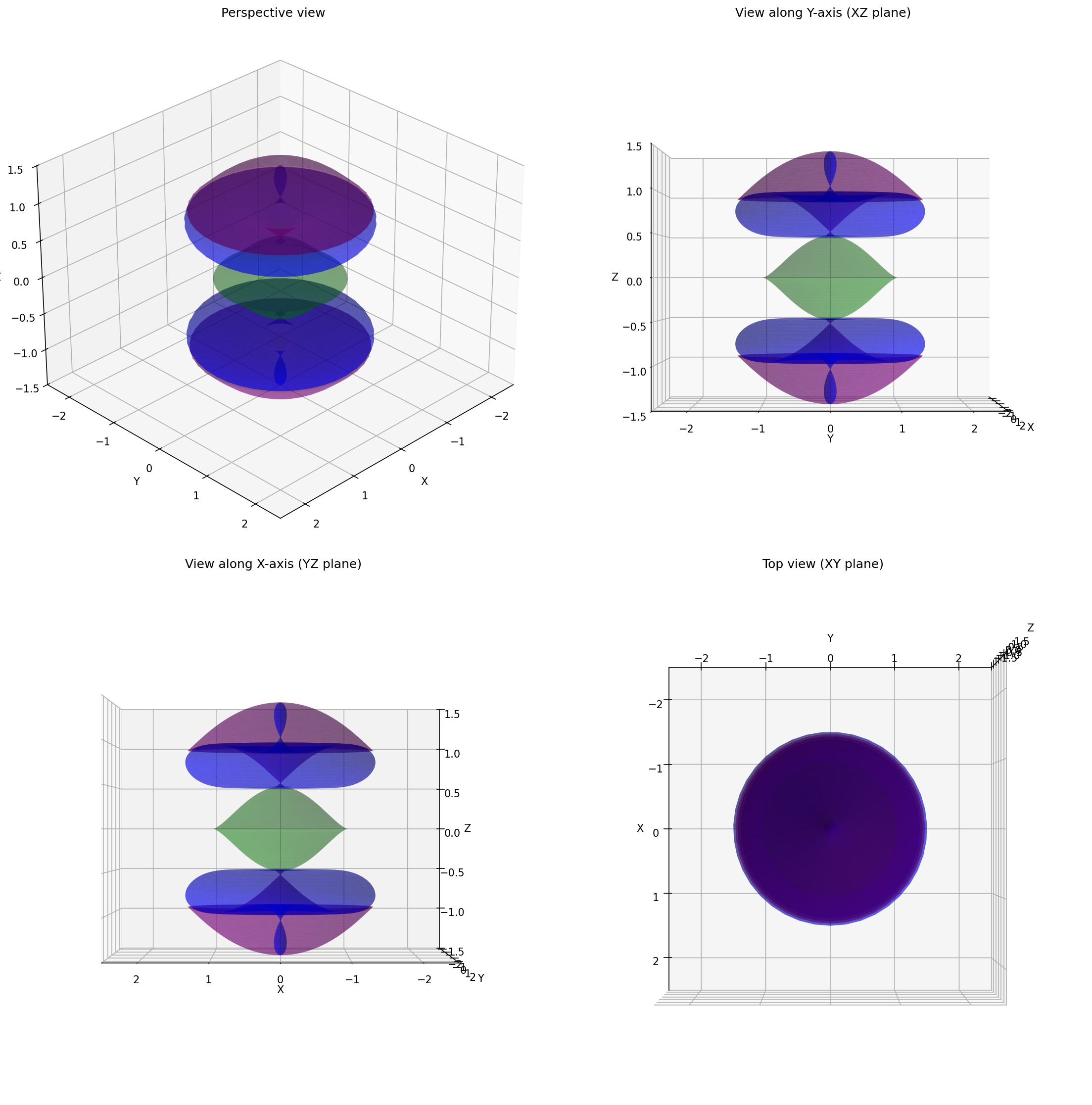}
 \caption{Implicit Plot} 
\end{subfigure}
\caption{Comparison of envelope surfaces: parametric representation (equations 9 and 10) versus implicit Cartesian form obtained via Sylvester resultant elimination. The parametric form traces specific envelope curves, while the implicit polynomial equations capture the complete algebraic variety, including additional solution branches arising from the elimination process.}
\label{fig:parametric_implicit_comparison}
\end{figure}

\textbf{Algebraic Verification.}
Algebraic verification was performed by substituting the parametric expressions $\rho^2(u) = \cos^2(u)(1 - 3\sin(u)/\sqrt{R})^2$ and $z(u) = (1/2)\sin^3(u) \pm 2/\sqrt{R}$ into the implicit equations. Numerical evaluation over 1000 parameter values confirmed $|F(\rho^2(u), z(u))| < 10^{-7}$, demonstrating equivalence at machine precision. Graphical comparison (Figure~\ref{fig:parametric_implicit_comparison}) further confirms geometric agreement. Notably, both envelope surfaces share identical horizontal projections in the $xy$-plane, differing only in their vertical displacement $\pm 2/\sqrt{R}$ from the clam curve.

\begin{figure}[htb]
\centering
\begin{subfigure}{5cm}
 \includegraphics[width=5cm]{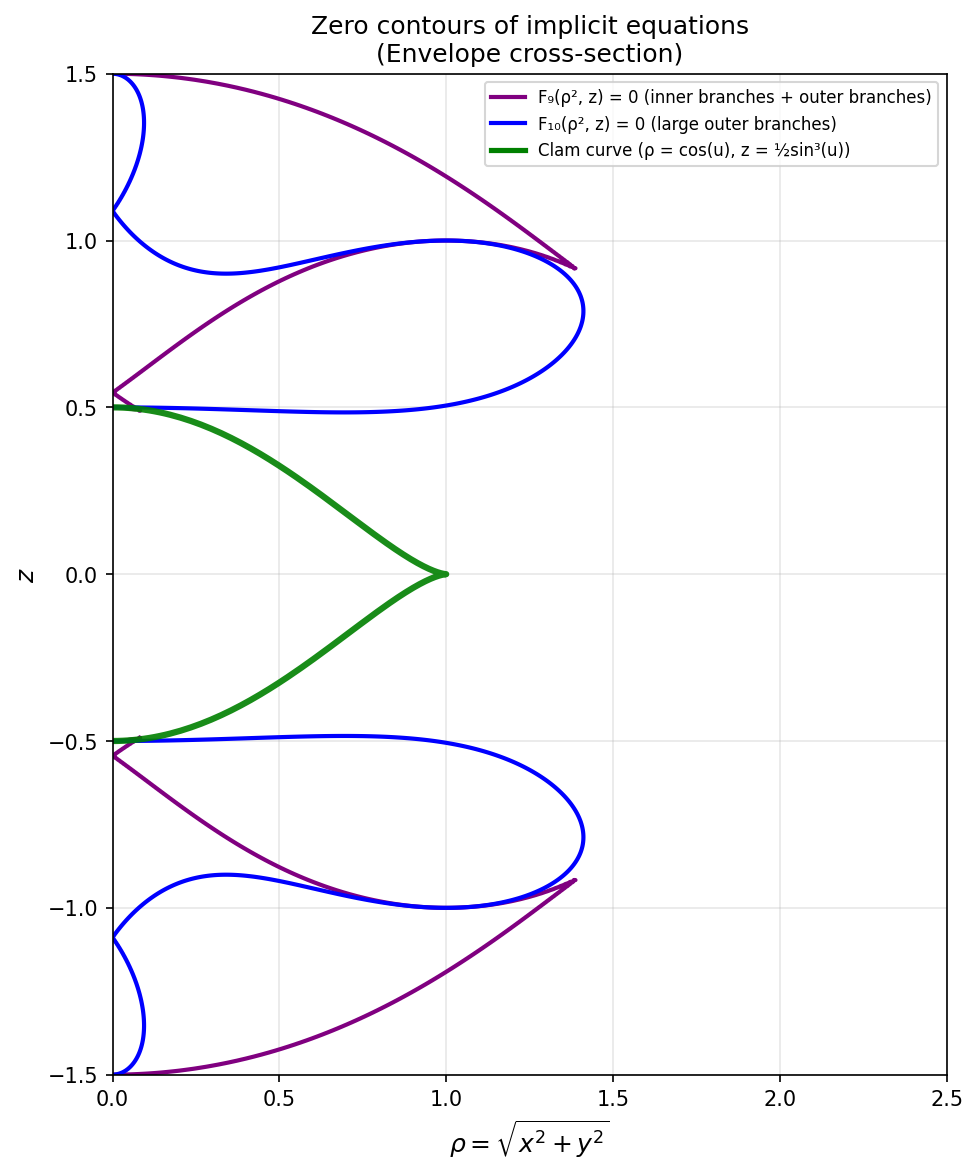}
\caption{} 
\end{subfigure}
\qquad
  \begin{subfigure}{5cm}
  \includegraphics[width=5.5cm]{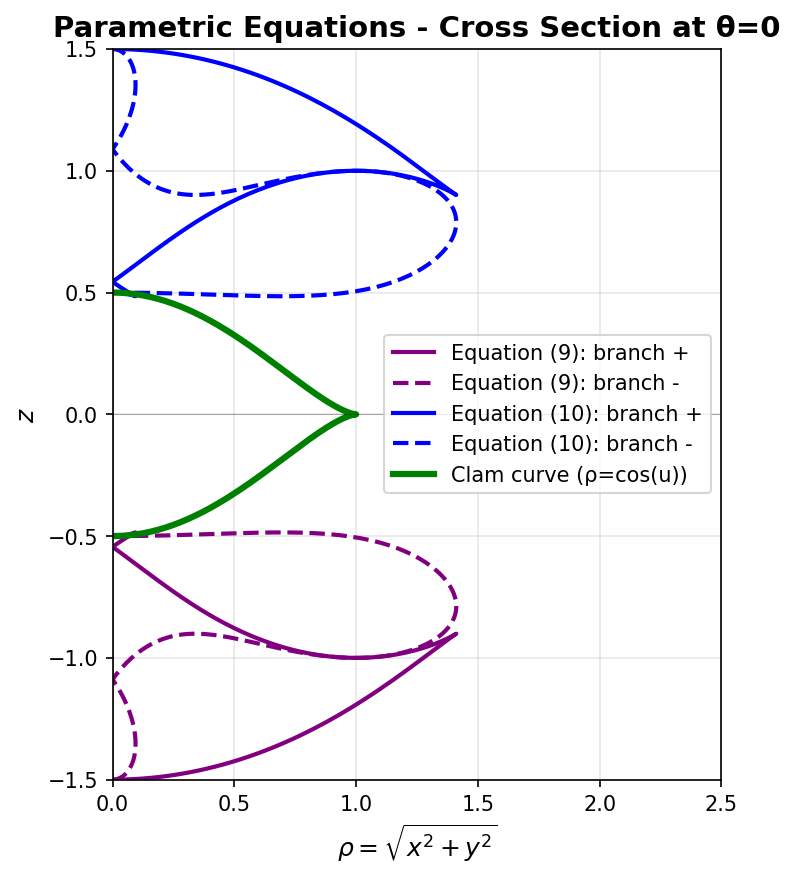}
 \caption{} 
\end{subfigure}
\caption{Plane sections of the clam and the offsets. This provides a visualization of how the envelope connects to the clam. }
\label{2D-section of the offsets}
\end{figure}

Greater insight is provided by considering a 2D section of the surfaces. This is displayed in Figure \ref{2D-section of the offsets}.

\section{Transition to 3D: offsets of the clam}
\label{section transition to 3D}
Suppose that a surface $\mathcal{S}$ is given by a cartesian equation $F(x,y,z)=0$. It is non singular at a point  $A$ if the gradient vector $\nabla F(x,y,z) \neq \vec{0}$ at $A$. A unit normal vector to $\mathcal{S}$ at $A$ is $\vec{N_A}= \frac {1}{\| \nabla F(x,y,z) \|} \nabla F(x,y,z)$.

We now generalize the definition of an offset to surfaces in 3D space.

\begin{definition}
\label{def offset surface}
Let $\mathcal{S}$ be a surface in the 3-dimensional space and $d$ a positive real number. At any point $A \in \mathcal{S}$ where the surface is smooth, we denote by $\vec{N_A}$ a unit normal vector. The points $B$ and $C$ are defined by $\overrightarrow{AB}=d \vec{N_A}$ and  $\overrightarrow{AC}=-d \vec{N_A}$. The union of the geometric loci of $B$ and $C$ when $A$ varies is called the offset at distance $d$ from the surface $\mathcal{S}$. One component is called the inner offset, while the other component is called the outer offset.
\end{definition}
It is possible, perhaps suitable, to define inner and outer offsets as we did in 2D. In fact, for large $d$, this may be counterintuitive.
 
\begin{remark}
A normal vector can also be computed as follows. Using a CAS, both methods yield the same answer. The user can choose which kind of commands to run. 
Suppose that the surface $\mathcal{S}$ (the progenitor) is defined by the parametrization $(x,y,z)=(\mathcal{S}_x(t), \mathcal{S}_y(t), \mathcal{S}_z(t) )$. A  normal vector $\vec{N}$ to the progenitor is given by $\vec{N} = \left( \frac{\partial x}{\partial v}, \frac{\partial y}{\partial v}, \frac{\partial z}{\partial v} \right) \times \left( \frac{\partial x}{\partial u}, \frac{\partial y}{\partial u}, \frac{\partial z}{\partial u} \right)$.
\end{remark}

If $d$ is the offset distance, the inner and outer offsets are determined respectively by the equations 
\footnotesize
\begin{equation}
\label{clam_offset_d}
\text{Offset}_{\text{in}}: \quad
\begin{cases}
x_{in} = \mathcal{S}_{\text{x}} + d \cdot {N}_{\text{x}}\\
y_{in} = \mathcal{S}_{\text{y}} + d \cdot {N}_{\text{y}}\\
z_{in} = \mathcal{S}_{\text{z}} + d \cdot {N}_{\text{z}}
\end{cases}
\quad \text{and} \quad  
\text{Offset}_{\text{out}}: \quad
\begin{cases}
x_{out} = \mathcal{S}_{\text{x}} - d \cdot {N}_{\text{x}}\\
y_{out} = \mathcal{S}_{\text{y}} - d \cdot {N}_{\text{y}}\\
z_{out} = \mathcal{S}_{\text{z}} - d \cdot {N}_{\text{z}}
\end{cases}
\end{equation}
\normalsize

Denoting $R_1=\sqrt{9\sin^4 u \cos^4 u \cos^2 v  + 9\sin^4 u \cos^4 u \sin^2 v + \sin^2 2u }$ and after simplification, Equations (\ref{param clam_offset_2a}) and (\ref{param clam_offset_2b}) read as follows:
\footnotesize
\begin{equation}
\label{param clam_offset_2a}
\text{Offset}_{\text{out}}
\begin{cases}
x_{out} = \cos v \cdot \left(\frac{3 d \sin^2 u \cos^2 u}{R_1} + \cos u \right)\\
y_{out}= \sin v  \cdot \left(\frac{3 d \sin^2 u \cos^2 u}{R_1} + \cos u \right)\\
z_{out} = \frac{d \sin 2u}{R_1} + \frac{\sin^3 u}{2}
\end{cases}
\end{equation}
\normalsize
and
\footnotesize
\begin{equation}
\label{param clam_offset_2b}
\text{Offset}_{\text{in}} =
\begin{cases}
x_{in} = \cos v \cdot \left(-\frac{3 d \sin^2 u \cos^2 u}{R_1} + \cos u \right)\\
y_{in} = \sin v  \cdot \left(-\frac{3 d \sin^2 u \cos^2 u}{R_1} + \cos u \right)\\
z_{in} = -\frac{d \sin 2u}{R_1} + \frac{\sin^3 u}{2}
\end{cases}
\end{equation}
\normalsize

 \begin{figure}[htb]
\centering
\begin{subfigure}{3.5cm}
 \includegraphics[width=3.5cm]{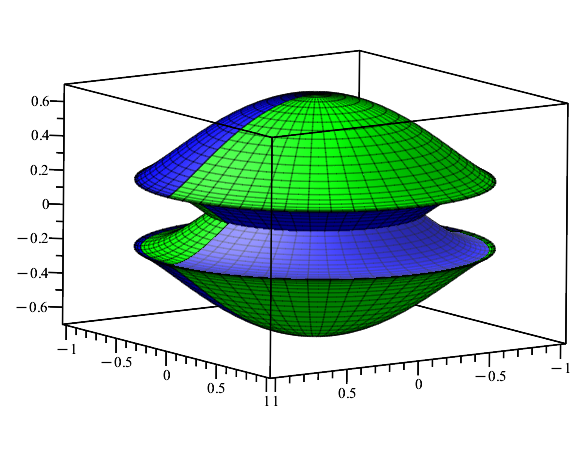}
\caption{$d = 0.2$} 
\label{offset_0_2}
\end{subfigure}
\qquad
  \begin{subfigure}{3.5cm}
  \includegraphics[width=3.5cm]{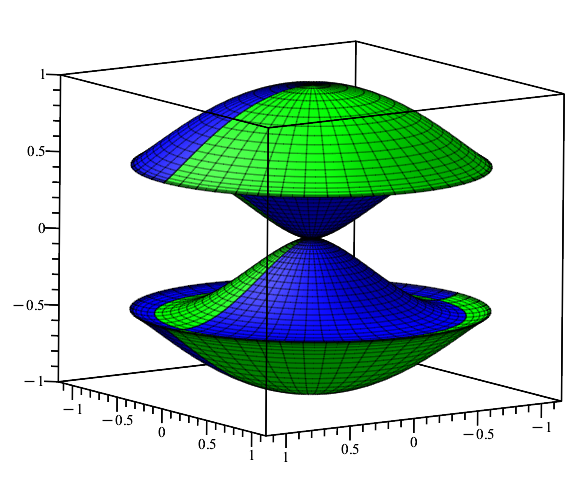}
 \caption{$d = 0.5$} 
\label{offset_0_5}
\end{subfigure}
\qquad
\begin{subfigure}{3.5cm}
 \includegraphics[width=3.5cm]{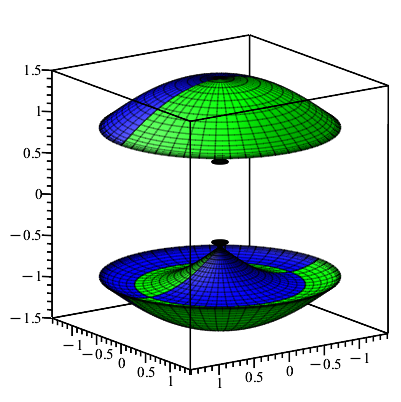}
 \caption{$d = 1$}
\label{offset_1}
 \end{subfigure}
\qquad
\begin{subfigure}{3cm}
 \includegraphics[width=3.5cm]{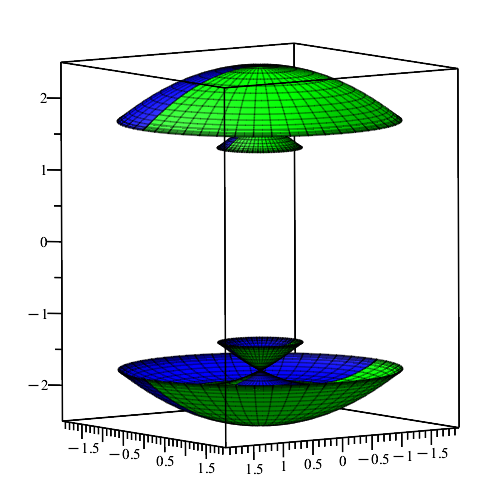}
\caption{$d = 2$} 
\label{offset_2}
\end{subfigure}
\qquad
\begin{subfigure}{3cm}
 \includegraphics[width=3.5cm]{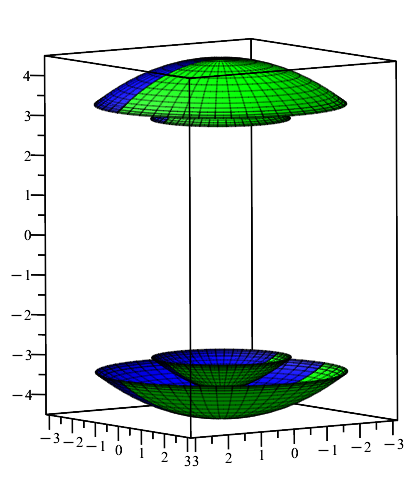}
\caption{$d = 4$} 
\label{offset_4}
\end{subfigure}
 \caption{Offsets of a clam at different distances}
\label{fig offsets}
\end{figure}

As for the kiss curve (and of course for numerous other occurences), plotting the equations of the offsets requires fine tuning with the choice of the parameters' intervals. For example, plotting Equations (\ref{param clam_offset_2a}) and (\ref{param clam_offset_2b}) for $u,v \in [0,2\pi ]$  yields irrelevant parts. Offsets after the fine tuning are displayed in Figure \ref{fig offsets} for several distances.

\subsection{Singularities of the offsets}
The graphical exploration reveals that in the 3D setting, the singularities are not isolated anymore. They are characterized by the condition:
$N = 0$. Solving the corresponding equations for $u$ and $v$ is a heavy task. For the present example, with its numerous symmetries, a different method can be applied: the surface under study can be obtained by revolving a kiss curve in the $xz-$plane by  about the $z-$axis. The singular points are also rotated, yielding a circle for each one. If a singular point of the kiss curve is  $P = (p_x,p_z)$, then revolving it  about the $z-$axis, we obtain 
\begin{equation}
\mathcal{C} = \begin{pmatrix} \cos(\theta)&-\sin(\theta)&0\\\sin(\theta)&\cos(\theta)&0\\0&0&1\end{pmatrix} \cdot \begin{pmatrix}p_x\\0\\p_z\end{pmatrix}
\end{equation}

Thus, for each singular point on the kiss curve,  the corresponding curve on the clam's offset is the geometric locus of points whose coordinates are : 
\begin{equation}
\begin{cases}
x = p_x \cdot \cos (\theta)\\
y = p_x \cdot \sin (\theta)\\
z = p_z
\end{cases}
\end{equation}
when $\theta \in [0,2 \pi ]$, i.e., circles. An example of these circles, both for the cusps and crunodes of the kiss curve, is shown in Figures \ref{3D_cuspidal_edges} and \ref{3D_crunodes_curves} on the offset surfaces.  If a crunode of the kiss curve lies on the axis of revolution, it remains invariant, i.e. an isolated point. As we have found the cusps and crunodes with 2 different methods for the kiss curves offset, we obtain the cuspidal edges and crunode curves separately. For example, for $d = 2$ we obtain the following singularities:
\begin{figure}[htb]
\begin{center}
\includegraphics[width=3cm]{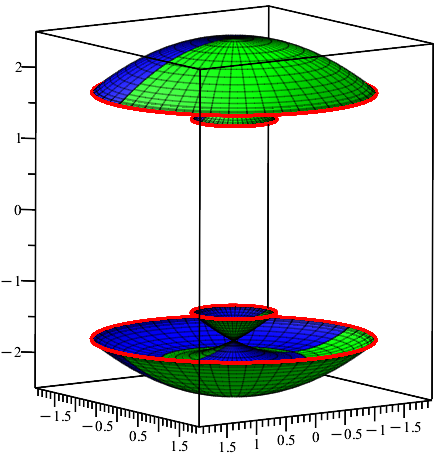}
\qquad
\includegraphics[width=3cm]{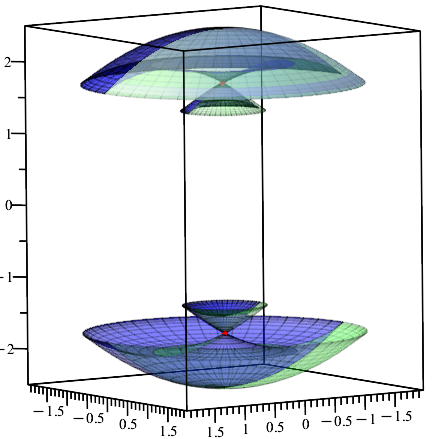}
\caption{Singular curves for virtual clam offsets when $d = 2$}
\label{3D_cuspidal_edges}
\end{center}
\end{figure}

\begin{figure}[htb]
\begin{center}
\includegraphics[width=3cm]{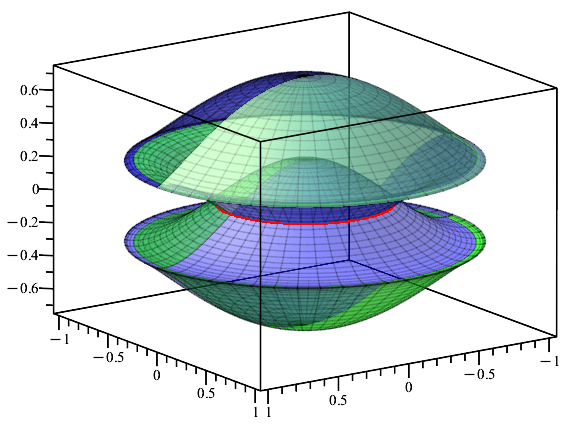}
\caption{Singular curves for virtual clam offsets when $d = 0.2$}
\label{3D_crunodes_curves}
\end{center}
\end{figure}

In order to find the singular curves, such as in Figures \ref{3D_cuspidal_edges} and \ref{3D_crunodes_curves}, we first find the 2D singular points of the kiss curves offsets (details are given in \cite{kiss curve}), and revolve them about the $z-$axis. We use the following code:

\footnotesize
\begin{verbatim}
eps := 1e-7;
singularities_plots := [];
for i to nops(D_singular_points) do
    pnt := D_singular_points[i];
    r := pnt[1];
    z := pnt[2];
    if abs(r) < eps then
        p := plot3d([r, r, z], scaling = constrained, color = red, 
        style = point, thickness = 20);
        singularities_plots := [op(singularities_plots), p];
    else
        p := plot3d([r*cos(t), r*sin(t), z], t = 0 .. 2*Pi, 
        scaling = constrained, color = red,
              style = line, thickness = 20);
        singularities_plots := [op(singularities_plots), p];
    end if;
end do;
display(singularities_plots);
\end{verbatim}
\normalsize
Note that, in the 2D setting, the kiss curve and its offsets are symmetric about the $y-$axis, so are the singular points. For these pairs of points, assuming that $P_1 = (p_x,p_z)$ and $P_2 = (-p_x,p_z)$, we obtain the following parameterizations for the cuspidal edges in the 3D setting:
\begin{equation}
\begin{cases}
x = p_x \cdot \cos (\theta)\\
y = p_x \cdot \sin (\theta)\\
z = p_z
\end{cases}
and
\begin{cases}
x = -p_x \cdot \cos (\theta)\\
y = -p_x \cdot \sin (\theta)\\
z = p_z
\end{cases}
\end{equation}
Both parameterizations describe the same curve. Thus, if the kiss curve offset has $\mathcal{N}$ singular points that are not located on the z-axis, the clams offset has $\frac{\mathcal{N}}{2}$ circles as singular curves.

\section{Final remarks}

\subsection{Using AI}
\label{using AI}
We tested several generative AI systems for computing envelopes and offsets. The answers were generally unreliable, as described in \cite{AI-open,AI-ggb,emprin-richard,AI-DPH}. However, we successfully employed one system after refining it through numerous iterations until it produced reliable results. As mentioned above, this process will be explained in detail in a subsequent paper.

\subsection{CAS assited determination of singular points}
If a surface $\mathcal{S}$ is given by a polynomial equation $F(x,y,z)=0$, then the singular points are determined by the solutions of the system of equations $\frac{\partial F}{\partial x}= \frac{\partial F}{\partial y}=\frac{\partial F}{\partial z}=0$. This system is generally solved easily using a CAS. Classifying the singularities is not trivial. This can be done using a parametrization of the surface, also a non trivial task. Maple has a built-in command for this. Then successive derivatives have to be computed and checked for vanishing or for linear independence sometimes. The computations are heavy and unilluminating.  

\subsection{Critical transition from 2D to 3D}
This approach for the study of singular points of plane curves obtained by various constructions has been widely used; for example,  \cite{estrella,kiss curve}. The transition from 2D to 3D can be quite easy, as in \cite{transition DPZ}. Here, as shown in section \ref{section transition to 3D} and subsection \ref{using AI}, the transition from 2D geometry  to 3D geometry induced a need for a transition between kinds of software. This transition is a critical one, in the sense of \cite{yerushalmy}. If in the past, networking between DGS and CAS was enough, as for example in \cite{safety,kiss curve,pedals}, the algebraic computations for an offset of a clam requested mastering a new technology.

\end{document}